\documentclass[12pt,a4paper]{article}
\usepackage{amssymb}

\usepackage{graphicx,amssymb,amsfonts,epsfig,amsthm,a4,amsmath,url}
\usepackage[latin1]{inputenc}

\newtheorem{ThmIntro}{Theorem}

\newtheorem{thm}{Theorem}[section]
\newtheorem{cor}[thm]{Corollary}
\newtheorem{lem}[thm]{Lemma}

\newtheorem{prop}[thm]{Proposition}
\theoremstyle{definition}
\newtheorem{defn}[thm]{Definition}

\theoremstyle{remark}

\numberwithin{equation}{section}

\newcommand{\Z}{\mathbf{Z}}

\newcommand{\N}{\mathbf{N}}
\newcommand{\R}{\mathbf{R}}

\newcommand{\Sl}{\text{SL}}
\newcommand{\Diam}{\text{Diam}}
\newcommand{\supp}{\text{Supp}}

\newcommand{\bpr}{\noindent \textbf{Proof}: ~}

\newcommand{\epr}{~$\blacksquare$}

\title{On the $L^p$-distortion of finite quotients of amenable groups.}
\author{Romain Tessera\footnote{This work was conducted while the author was visiting
the Bernoulli center in Lausanne. The author is supported by the
NSF grant DMS-0706486.} }
\date{\today}

\begin{document}

\baselineskip=16pt

\maketitle

\begin{abstract}
We study the $L^p$-distortion of finite quotients of amenable
groups. In particular, for every $2\leq p<\infty$, we prove that
the $\ell^p$-distortions of the groups $C_2\wr C_n$ and
$C_{p^n}\ltimes C_n$ are in $\Theta((\log n)^{1/p}),$ and that the
$\ell^p$-distortion of $C_n^2\ltimes_A \Z$, where $A$ is the
matrix $\left(\begin{array}{cc}
2 & 1 \\
1 & 1
\end{array}\right)$ is in $\Theta((\log \log n)^{1/p}).$
\end{abstract}



\section{The main results}

\subsection{Distortion}

Let us first recall some basic definitions.

\begin{defn}

\

\begin{itemize}
\item Let $0<R\leq \infty$. The distortion at scale $\leq R$ of an
injection between two discrete metric spaces $F:(X,d)\to (Z,d)$ is
the number (possibly infinite)
$$dist_R(F)=\sup_{0<d(x,y)\leq R}\frac{d(f(x),f(y))}{d(x,y)}\cdot\sup_{0<d(x,y)\leq R}\frac{d(x,y)}{d(f(x),f(y))}.$$
If $R=\infty$, we just denote $dist(F)$ and call it the distortion
of $F$.

\item The $\ell^p$-distortion $c_p(X)$ of a finite metric space
$X$ is the infimum of all $dist_F$ over all possible injections
$F$ from $X$ to $\ell^p$.
\end{itemize}
\end{defn}

Let $G$ be a finitely generated group. Let $S$ be a symmetric
finite generating subset of $G$. We equip $G$ with the
left-invariant word metric associated to $S$:
$d_S(g,h)=|g^{-1}h|_S=\min\{n\in \N, g^{-1}h\in S^n\}.$ Let
$(G,S)$ denote the associated Cayley graph of $G$: the set of
vertices is $G$ and two vertices $g$ and $h$ are joined by an edge
if there is $s\in S$ such that $g=hs$. Note that the graph metric
on the set of vertices on $(G,S)$ coincides with the word metric
$d_S$.

Let $\lambda_{G,p}$ denote the regular representation of $G$ on
$\ell^p(G)$ for every $1\leq p\leq \infty$ (i.e.
$\lambda(g)f(x)=f(g^{-1}x)$). The $\ell^p$-direct sum of $n$
copies of $\lambda_{G,p}$ will be denoted by $n\lambda_{G,p}$.

Our main results are the following theorems.

\begin{ThmIntro}\label{lamplighterThmIntro}
Let $m$ be an integer $\geq 2$. For all $n\in \N$, consider the
finite lamplighter group $C_m\wr C_n=(C_m)^{C_n}\ltimes C_n$
equipped with the generating set $S=((\pm 1_{0},0), (0,\pm 1))$,
where $1_0\in (C_m)^{C_n}$ is the characteristic function of the
singleton $\{0\}$. For every $2\leq p< \infty$, there exists
$C=C(p,m)<\infty$ such that
$$C^{-1}(\log
n)^{1/p} \leq c_p(C_2\wr C_n,S)\leq C(\log n)^{1/p}.$$
\end{ThmIntro}
Note that the upper bound has been very recently proved for $p=2$
by Austin, Naor, and Valette \cite{ANV}, using representation
theory. The proof that we propose here is shorter and
completely elementary. On the other hand, the lower bound was
known (see \cite{LNP}, or Section~\ref{UpperboundSection}).

\begin{ThmIntro}\label{BSThmIntro}
Let $m$ be an integer $\geq 2$. For all $n\in \N$, consider the
group $BS_{m,n}=C_{m^n}\ltimes C_n$ equipped with the generating
set $S=\{(\pm 1,0), (0,\pm 1)\}$. For every $2\leq p< \infty$,
there exists $C=C(p,m)<\infty$ such that
$$C^{-1}(\log
n)^{1/p} \leq c_p(G_n,S)\leq C(\log n)^{1/p}.$$
\end{ThmIntro}

\begin{ThmIntro}\label{SolThmIntro}
For all $n\in \N$, consider the group $SOL_n=C_{n}\ltimes_A
C_{o(A,n)}$, where $A$ is a matrix of $\Sl_2(\Z)$ with eigenvalues
of modulus different from $1$, e.g. the matrix
$\left(\begin{array}{cc}
2 & 1 \\
1 & 1
\end{array}\right)$, and where $o(A,n)$ denotes the order of $A$ in
$SL_2(C_n)$. Equip $G$ with the generating set $S=\{(\pm 1,0),
(0,\pm 1)\}$. For every $2\leq p<\infty$, there exists
$C=C(p)<\infty$ such that
$$C^{-1}(\log \log n)^{1/p} \leq c_p(G_n,S)\leq C(\log \log n)^{1/p}.$$
\end{ThmIntro}

\subsection{About the constructions}\label{about-construction-Section}
We will say that map $F: G\to E$ from a group $G$ to a Banach
space is equivariant if it is the orbit of $0$ of an isometric
affine action of $G$ on $E$. Let $\sigma$ be such an action. The
equivariance of  $F(g)=\sigma(g).0$ implies that
$\|F(g)-F(h)\|=\|F(g^{-1}h)\|$. Hence the distortion at scale
$\leq R$ of $F$ is just given by
$$dist_R(F)=\sup_{0<|g|_S\leq R}\frac{|g|_S}{\|F(g)\|}\cdot\sup_{0<|g|_S\leq R}\frac{\|F(g)\|}{|g|_S}.$$

All the groups involved in the main theorems are of the form
$G=N\ltimes A$ where $A$ is a finite cyclic group. To prove an
upper bound on $c_p(G),$ our general approach is to construct an
embedding $F=F_1\oplus^{\ell^p}F_2$, where $F_1$ is the orbit of
$0$ of an affine action $\sigma_1$ of $G$, whose linear part is
$K\lambda_{G,p}$ (for some $K\in \N$), and such that for $R=
\Diam(N)$, we have $$dist_R(F_1)\approx (\log R)^{1/p}.$$ More
precisely, for $F_{m,n}$ and $BS_{m,n}$ (resp. for $SOL_{A,n}$),
we will need $K\approx \log(mn)$  (resp. $K\approx \log\log n$)
copies of $\lambda_{G,p}$.

For $G=F_{m,n}$ or $BS_{m,n}$, we can take $F=F_1$ since
$\Diam(N)\approx \Diam(G)\approx n$ (see
Proposition~\ref{propDiam}). But, for $G=SOL_{A,n}$, we have
$\Diam(N)\approx \log n$, which can be much less than
$\Diam(G)\approx o(A,n)$. Hence, the solution in this case is to
add some map $F_2: G/N\approx C_{o(A,n)}\to \ell^p$ with a bounded
distortion (for instance, take the orbit of $0$ under the action
of $C_{o(A,n)}$ on $\R^2$ such that $1$ acts by rotation of center
$(o(A,n),0)$ and angle $2\pi/o(A,n)$).

Note that Theorem~\ref{SolThmIntro} also holds for the group $C_{n}\ltimes_A
\Z$, in which case we can take an action of $\Z$ by translations on $\R$ to embed the quotient with bounded distortion (i.e. for $F_2$).
\section{Upper bounds on the distortion}\label{UpperboundSection}

Let $1\leq p\leq \infty.$ Recall \cite{T1} that the
left-$\ell^p$-isoperimetric profile in balls of $(G,S)$ is defined
by
$$J_{G,S,p}(n)=\sup_{\supp(f)\subset B(1,n)}\frac{\|f\|_p}{\sup_{s\in S}\|\lambda(s)f-f\|_p},$$
where $B(1,n)$ denotes the open ball of radius $n$ and center $1$
in $(G,S)$. For convenience, we will

Our main result in \cite{T1} consisted in showing that a lower
bound on the isoperimetric profile can be used to construct
metrically proper affine isometric actions of $G$ on $\ell^p(G)$
whose compressions satisfy lower bounds which are optimal in certain cases.
Here, we will use it to produce upper bounds on the
$\ell^p$-distortion of finite groups.

On the other hand, as explained in \cite{T2}, if $X=(G,d_S)$ is a
Cayley graph, then the inequality $J_{p,G}\geq J$ for some non-decreasing
function $J:\R_+\to \R_+$ implies Property A(J,p) (see
\cite[Definition~4.1]{T2}) for the space $X$ (if the group $G$ is
amenable, a standard average argument actually shows that this is
an equivalence). So in a large extend, the results of the present
paper are easy consequences of the method explained in \cite{T2}.

A crucial remark is that $J_{G,S,p}$ is a local quantity, and
hence behaves well under quotients. Namely, we recall the
following easy fact.

\begin{prop}\label{easyprop} (for a proof, see
\cite[Theorem~4.2]{T3}) Let $\pi:G\to Q$ be a surjective
homomorphism between two finitely generated groups and let $S$ be
a symmetric generating subset of $G$. Then
$$J_{G,S,p}\leq J_{Q,\pi(S),p}.$$
\end{prop}

Our main technical tool is the following proposition, which is an
analogue of \cite[Proposition~4.5]{T2}. For the convenience of the
reader, we give its relatively short proof in
Section~\ref{sectionProof}.

\begin{ThmIntro}\label{mainThmIntro}
Let $X=(G,S)$ be a finite Cayley graph such that $J_{G,S,p}(r)\geq
J(r)$ when $r\leq R$, for some $R\leq \Diam(G)/2$. Then, there
exists an affine isometric action $\sigma$ of $G$ on such that
\begin{itemize}
\item the linear part of $\sigma$ is the $\ell^p$-direct sum of
$K=[\log R]$ regular representations of $G$ in $\ell^p(G)$.

\item The orbit of $0$ induces an injection $F: G\to
\bigoplus_{k=0}^{K-1}\ell^p(G)$ such that $$dist_R(F)\leq
2\left(2\int_2^{R/2}\left(\frac{t}{J(t)}\right)^p\frac{dt}{t}\right)^{1/p}.$$
In particular, if $J(t)= t/C$, then
$$dist_R(F)\leq 2C\left(2\log (R/2)\right)^{1/p}.$$
\end{itemize}
\end{ThmIntro}

\begin{cor}\label{corIntro}
Assume that $G_n$ has diameter $\leq n$ and that $J_{G,p}(t)\geq
t/C$, then, $c_p(G_n)\leq 2C\left(2\log (n/4))\right)^{1/p}.$
\end{cor}

On the other hand, we have proved in \cite{T1} that the following
finitely generated groups satisfy $J_p(t)\geq t/C$ for some
$C<\infty$ and for all $1\leq p<\infty$.
\begin{itemize}
\item the lamplighter group $L_{m}=C_m\wr \Z$;

\item solvable Baumslag-Solitar groups $BS_m=\Z[1/m]\ltimes \Z$
for all $m\in \N$, where $n\in \Z$ acts by multiplication by
$m^n$;

\item polycyclic groups. Here, we will focus on the following
example: $SOL_A=\Z^2\ltimes_A\Z$ where $A$ is a matrix of
$\Sl_2(\Z)$ with eigenvalues of modulus different from $1$, e.g.
the matrix $\left(\begin{array}{cc}
2 & 1 \\
1 & 1
\end{array}\right).$
\end{itemize}
Note that respectively $L_{m,n}$, $BS_{m,n}$ and $SOL_{A,n}$ are
quotients of $L_m$, $BS_m$ and $SOL_A$.

\section{Proofs of the main theorems}
\subsection{Upper bounds}
Thanks to Corollary~\ref{corIntro}, the upper bounds in Theorems~\ref{lamplighterThmIntro}, \ref{BSThmIntro} and \ref{SolThmIntro} follow from the following upper bounds on the diameters of the groups
$L_{m,n}$, $BS_{m,n}$ and $SOL_{A,n}$ (for the latter, see the discussion in Section~\ref{about-construction-Section}).

\begin{prop}\label{propDiam}
We have
\begin{itemize}
\item[(i)] $\Diam(L_{m,n})\leq (m+3)n$;

\item[(ii)] $\Diam(BS_{m,n})\leq (m+1)n$;

\item[(iii)] Let $N_n\simeq C_n^2$ be the kernel of $SOL_{A,n}\to
C_{o(A,n)}$. Then, with the distance on $N_n$ induced by the word
distance on $SOL_{A,n}$, we have $\Diam(N_n)\leq c\log n$ for some
$c=c(A)>0$.
\end{itemize}
\end{prop}
\bpr For (i), see \cite{Pa}. For (ii), note that every element of
$C_{m^n}$ can be written as
$$\sum_{i=0}^{n-1}a_im^i=a_0+m(a_1+m(a_2+\ldots)\ldots),$$
where $0\leq a_i\leq m-1$. Finally, (iii) follows from the
following well known lemma.\epr

\begin{lem}\label{lemmaSOL}
Let $N\sim \Z^2$ be the kernel of $SOL_{A}\to \Z$. For all $r\geq
1$, denote by $B_{N,SOL_A}(r)$ (resp. $B_N(r)$), the ball of
radius $r$ for the metric on $N$ induced by the word length on
$SOL_A$ (resp. for the usual metric on $\Z^2$). There exists some
$\alpha=\alpha(A)<\infty$ such that
$$B_{N}(1,e^{r/\alpha})\subset B_{N,SOL_A}(r)\leq B_{N}(1,e^{\alpha r}).$$
\end{lem}
\bpr Note that $SOL_A$ embeds as a co-compact lattice in the
connected solvable Lie group $G=\R^2\ltimes_A R$, such that $N$
maps on a (co-compact) lattice of $\tilde{N}=\R^2$. The lemma
follows from the fact that $\tilde{N}$ is the exponential radical
of $G$ (Guivarc'h \cite{Gu} was the first one to introduce and to
study the exponential radical of a connected solvable Lie group, without
actually naming it, and this was rediscovered by Osin
\cite{Osin}). \epr

\subsection{Lower bounds}
To obtain the lower bound on the distortion, we will need the
following notion of relative girth.

\begin{defn}
Let $\pi:G\to Q$ be a surjective homomorphism between two finitely
generated groups and let $S$ be a symmetric generating subset of
$G$. Denote by $X=(G,S)$ and $Y=(H,\pi(S))$. The relative girth
$g(Y,X)$ of $Y$ with respect to $X$ is the maximum integer $n\in
\N$ such that a ball of radius $n$ in $Y$ is isometric to a ball
of radius $n$ in $X$.
\end{defn}

Recall \cite{Bourgain} that the rooted binary tree $T_n$ of dept
$n$ satisfies $c_p(T_n)\geq c(\log n)^{1/p}$ for all $2\leq
p<\infty$ and for some constant $c>0$. The following remark
follows trivially from this result and from the definition of
relative girth.
\begin{prop}
We keep the notation of the previous definition. Assume that $X$
contains a bi-Lipschitz embedded 3-regular tree. Then there exists
some $c>0$ such that $c_p(Y)\geq c(\log g(X,Y))^{1/p}.$\epr
\end{prop}

On the other hand, the groups $L_m$, $BS_m$ and $SOL_A$ are
solvable non-virtually nilpotent. Hence by \cite{QI}, they admit a
bi-Lipschitz embedded 3-regular tree (for the lamplighter, see
also \cite{LPP}). So to prove the lower bounds of
Theorems~\ref{lamplighterThmIntro}, \ref{BSThmIntro} and
\ref{SolThmIntro}, we just need to find convenient lower bounds
for the relative girths, which is done by the following
proposition.

\begin{prop}
We have
\begin{itemize}
\item[(i)] $g(L_{m,n},L_m)\geq n$;

\item[(ii)] $g(BS_{m,n},BS_m)\geq n$;

\item[(iii)] $g(SOL_{A,n},SOL_A)\geq c\log n$ for some $c=c(A)>0$.
\end{itemize}
\end{prop}
\bpr The only non-trivial case, (iii), follows from
Lemma~\ref{lemmaSOL}. \epr

\section{Proof of Theorem~\ref{mainThmIntro}}\label{sectionProof}

Let $f_0$ be the dirac at $1$, and for every integer $1\leq k\leq
K$, choose a function $f_k\in \ell^p(G)$ such that
\begin{itemize}
\item the support of $f_k$ is contained in the ball $B(1,2^k)$,

\item $\|f_k\|_p\geq J(2^k)$

\item $\sup_{s\in S}\|\lambda(s)f_k-f_k\|_p\leq 1$
\end{itemize}
For all $v=(v_k)_{1\leq k\leq n}\in K\ell^p(G)$ and all $g\in G$,
define
$$\sigma(g)v=\bigoplus_k^{\ell^p}(\lambda(g)v_k+F_k)$$ where
$$F_k(g)=\left(\frac{2^k}{J(2^k)}\right)(f_k-\lambda(g)f_k).$$

Now consider the map $F=\bigoplus^{\ell^p}b_k: G\to K\ell^p(G)$.
For all $g\in G$, we have
\begin{eqnarray*}
\|F(g)\|_p &  = & \|b(g)\|_p \\  & \leq
& \left(\sum_{k=0}^n
\left(\frac{2^k}{J(2^k)}\right)^p\|\lambda(g)f_k-f_k|\|_p^p\right)^{1/p} \\
& \leq & \left(\sum_{k=0}^n
\left(\frac{2^k}{J(2^k)}\right)^p\right)^{1/p}\\
&\leq &
|g|_S\left(\int_1^{\Diam(G)/2}\left(\frac{t}{J(t/2)}\right)^p\frac{dt}{t}\right)^{1/p}\\
&=&
2^{2/p}|g|_S\left(\int_1^{\Diam(X)/4}\left(\frac{t}{J(t)}\right)^p\frac{dt}{t}\right)^{1/p}.
\end{eqnarray*}
On the other hand, since $f_k$ is supported in $B(1,2^k)$, if
$|g|_S\geq 2.2^k$, then the supports of $f_k$ and $\lambda(g)f_k$
are disjoint. Thus,
\begin{eqnarray*}
\|F(g)\|_p & =  & \|b(g)\|_p \\
 & \geq &
\|b_k\|_p \\
& =  & 2^{1/p}\frac{2^k}{J(2^k)}\|f_k\|_p
\\
& \geq  & 2^{1/p}2^k,
\end{eqnarray*}
whenever $d_S(x,y)\geq 2.2^k$. To conclude, we have to consider
the case when $g\in S\smallsetminus\{1\}$. But as $f_0$ is a dirac
at $1$, $\|F(g)\|_p\geq 1$. So we are done.\epr

\bigskip
\footnotesize

\noindent \noindent Romain Tessera\\
Department of mathematics, Vanderbilt University,\\ Stevenson
Center, Nashville, TN 37240 United,\\ E-mail:
\url{tessera@clipper.ens.fr}

\end{document}